# The Meaning of Infinity

W. Mückenheim

University of Applied Sciences, Baumgartnerstraße 16, D-86161 Augsburg, Germany

mueckenh@rz.fh-augsburg.de

**Abstract.** The notions of potential infinity (understood as expressing a direction) and actual infinity (expressing a quantity) are investigated. It is shown that the notion of actual infinity is inconsistent, because the set of all (finite) natural numbers which it is ascribed to, cannot contain an actually infinite number of elements. Further the basic inequality of transfinite set theory $\aleph_0 < 2^{\aleph_0}$ is found invalid, and, consequently, the set of real numbers is proven denumerable by enumerating it.

## 1. Introduction

Infinity has fascinated philosophers and mathematicians at all times. So it is not surprising that the ideas concerning this problem evolved quite differently, not only with respect to the fundamental question whether its meaning is a potential one or an actual one, but also with respect to the rules to be observed in the infinite. The fathers of calculus had to struggle with infinitely small quantities which, however, were not zero, so that bishop Berkeley ironically remarked: "They are neither finite quantities, nor quantities infinitely small, nor yet nothing. May we not call them the ghosts of departed quantities?" But more than the infinitely small, which could successfully be dealt with, at least in *applying* the calculus, the infinitely large resisted the human mind.

While, according to Galilei, the infinite should obey another arithmetic than the finite, Leibniz, on the contrary, expected the rules of finite mathematics to remain valid in the infinite. It seems reasonable, to a certain degree, to maintain approved rules because (1) this method always has turned out fruitful in history of mathematics when numbers and notions had to be generalised, and because (2) we have no *a priori* knowledge or imagination of the actual infinite, so we would be helpless without relying to the firm rules of finite mathematics. The principle of Leibniz is most distinctly realised in trusting in the unlimited validity of the law of complete induction as a basic law in classical mathematics as well as in non-standard analysis, just because there is no natural number which could not be member of a finite set. A question of similar importance is the meaning of a one-to-one correspondence or bijection between two sets. Albert of Saxonia proved already in the 14$^{th}$ century that a wooden bar of infinite length has the same contents as the whole world. Galilei proved the equivalence of the set of all natural numbers and the set of all squares. Both proofs were based on the assumption that two sets are equivalent or have the same cardinality if a bijection between them can be established. Bernard Bolzano, the first mathematician who made a substantial contribution to the science of the infinite and invented the notion of "set" (using the German expression "Menge"), denied this principle. Though agreeing that there are as many natural numbers as their squares, he asserted that a long line contains more points than a short line and that the set of focuses of all ellipses is twice as large as the set of their centres. Georg Cantor, on the contrary, supported and popularised the opinion that a bijection always, in the finite as well as in the infinite, proves equivalence. Convinced that actual infinity does really exist in three different realms, namely in God and his properties, in nature by infinitely many created entities, and in mathematics, he developed transfinite set theory based upon two central laws: (1) the possibility of a bijection between two sets proves their equivalence, and (2) the set of natural numbers $\mathbb{N}$ has an actually infinite number of elements, denoted by $\aleph_0$. He considers "the infinite set of all finite whole numbers"[1] as showing "the simplest example of an actually infinite quantity"[2], an opinion which had also been shared by Bolzano[3].

## 2. Potential infinity versus actual infinity

In order to analyse the notion of an actually infinite set of finite numbers in detail, we must sharply distinguish between potential and actual infinity. "The essential difference between the potential infinite and the actual infinite did, strange to say, not prevent several cases of intermingling both notions in the development of the newer mathematics", Cantor writes in a review of set theory[4] but his sentences just quoted[1,2] show that Cantor himself was not free of committing this error.

---

[1] Denn, daß es sich bei der unendlichen Menge aller endlichen ganzen Zahlen ... [1, p. 402]. Die aktual-unendliche Menge aller positiven, endlichen ganzen Zahlen ... [1, p. 412].

[2] ... wogegen die durch ein Gesetz begrifflich durchaus bestimmte Menge aller ganzen endlichen Zahlen das einfachste Beispiel eines aktual-unendlichen Quantums darbietet [1, p. 409].

[3] Wenn wir die Reihe der natürlichen Zahlen: 1, 2, 3, 4, 5, 6, ... betrachten: so werden wir gewahr, dass die Menge der Zahlen, die diese Reihe, anzufangen von der ersten (der Einheit) bis zu irgend einer z. B. der Zahl 6, enthält, immer durch diese letzte selbst ausgedrückt wird. Somit muss ja die Menge aller Zahlen genau so gross als die letzte derselben und somit selbst eine Zahl, also nicht unendlich sein. Das Täuschende dieses Schlusses verschwindet auf der Stelle, sobald man sich nur erinnert, dass in der Menge aller Zahlen in der natürlichen Reihe derselben keine die letzte stehe; dass somit der Begriff einer letzten (höchsten) Zahl ein gegenstandloser, weil einen Widerspruch in sich schliessender, Begriff sei [2, p. 21].

[4] Die wesentliche Verschiedenheit, welche hiernach zwischen den Begriffen des potentialen und aktualen Unendlichen besteht, hat es merkwürdigerweise nicht verhindert, daß in der Entwicklung der neueren Mathematik mehrfach Verwechslungen beider Ideen vorgekommen sind, derart, daß in Fällen, wo nur ein potentiales Unendliches vorliegt, fälschlich ein Aktual-Unendliches angenommen wird, oder daß umgekehrt Begriffe, welche nur vom Gesichtspunkte des aktualen Unendlichen einen Sinn haben, für ein potentiales Unendliches gehalten werden [1, p. 409].

If a set consists of different finite whole numbers only, then it has necessarily a finite number of members and, hence, is a finite set. Consider the sequence 1, 2, 3, ..., *n* with finite *n*. *n* is equal to the number of terms. This situation does not change with growing *n*, because *all* $n \in \mathbb{N}$ are finite. If we count 1, 2, 3, ..., then we are in the infinite, and further great expectations are unjustified. "Infinite" can here only be used in its potential meaning, namely of *n* "being capable of growing beyond any threshold" though always remaining a finite number. $\mathbb{N}$ cannot contain an *actually infinite* number of members, i.e., a number of members which is larger than *any* member of the set, because the natural numbers enumerate themselves. Position *n* and magnitude *f(n)* = *n* are identical, without any single exception.

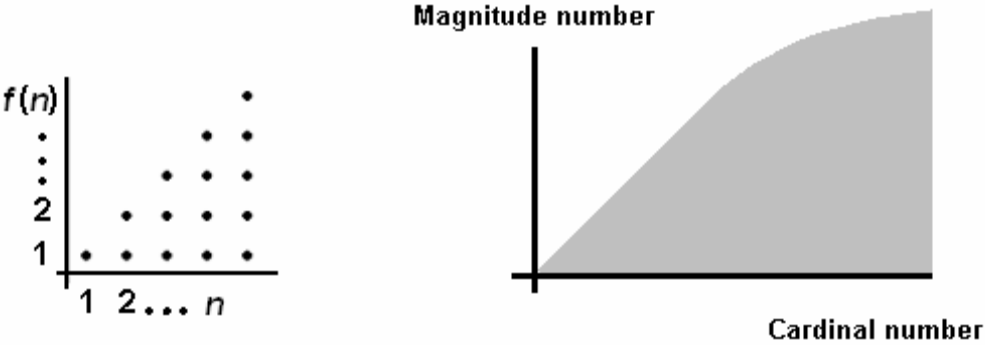

**Fig. 1.** The function *f(n) = n*.     **Fig. 2.** An infinite set of finite numbers.

And it is unreasonable to expect or to require that the linear function *f(n) = n*, shown in Fig. 1, "somewhere in the infinite" should change its slope as depicted in Fig. 2, or should exhibit a gap. Therefore it is inconsistent to speak of an infinite set of finite numbers in its actual meaning.

But in order to avoid a discussion about the *last number* of $\mathbb{N}$ which, as Bolzano remarks [2, p. 21], does not exist, take the following innocent

Theorem: *Every finite nonempty set of positive even numbers contains at least one number n surpassing the cardinal number of the set.*

It is easy to see that this property is shared by at least half of those even numbers, and it is an easy exercise to prove this theorem by complete induction. Therefore it holds for *all* sets of positive even numbers including the set of *all* positive even numbers, should it exist, because there is absolutely no one which could not be an element of a finite set. This proves the existence of a natural number $n \in \mathbb{N}$ surpassing the cardinal number of the set of all positive even numbers $\aleph_0$, although this smallest transfinite cardinal number, according to its inventor Cantor, should be larger than any finite number.

Our theorem could be sharpened by using prime numbers instead of even numbers but already from the original version of this simple proof we see clear enough that it is inconsistent to speak of *an infinite set of finite numbers*. Finite numbers can only form a potentially infinite set. An actually infinite set cannot exist other than including its cardinal number $\aleph_0$ or the number "ordinal infinity", denoted by $\omega$. This had been unconsciously acknowledged by Cantor himself already: "*Every number smaller than $\omega$ is a finite number, and its magnitude is surpassed by other finite numbers $\nu$*"[5]. Here the phrase "by other finite numbers" is obviously to be interpreted as "by such finite numbers which did not yet belong to the set containing *every finite number*". $\omega$ *is not a natural number*, because instead of the fundamental inequality obeyed by natural numbers, $n < n + 1$, it obeys the different law $\omega = 1 + \omega$. Therefore an actually infinite set of *only* natural numbers, in particular the complete set of *exclusively all* natural numbers, usually abbreviated by $\mathbb{N}$, cannot and does not exist at all (unless we choose one and the same number infinitely often, a method which Cantor certainly would not have approved), and the onus of finding a natural number $n$ which does not take its finite place $n$ in the sequence, rests on those who support this self-contradicting notion. Of course, we can denote this set by $\mathbb{N}$ as we can name dwarfs, fairies and unicorns. But that does not found its existence.

Cantor[6] himself admits that there are impossible or inconsistent sets. One of them is the set of all sets, because it must but cannot contain its own power set. The set of all natural numbers, if indeed actually existing, must but cannot contain the arithmetic sum or product of all of its elements. But there is no logical difference: Why should the elements of an actually existing infinite set be capable of generating their power set but not their arithmetic sum?

---

[5] Jede kleinere Zahl als $\omega$ ist eine endliche Zahl und wird von anderen endlichen Zahlen $\nu$ der Größe nach übertroffen [1, p. 406].
[6] Es gibt also bestimmte Vielheiten, die nicht zugleich Einheiten sind, d. h. solche Vielheiten, bei denen ein reales "Zusammensein aller ihrer Elemente" unmöglich ist. Diese sind es, welche ich "inkonsistente Systeme", die anderen aber "Mengen" nenne [1, p. 448].

## 3. Impossible bijections

As we have seen, actual infinity does not exist, but we can show by an independent proof that the next known cardinal number would cause a contradiction too. Cantor proved that the cardinal number $2^{\aleph_0}$ of the set $\mathbb{R}$ of all real numbers is strictly larger than $\aleph_0$ [1, pp. 115, pp. 139, pp. 278]. This can briefly be expressed by the inequality [1, p. 448]

$$\aleph_0 < 2^{\aleph_0}. \tag{1}$$

But it is easy to see that mathematics sometimes requires the existence of not less than $2^{\aleph_0}$ natural numbers. Bishop Nicole de Oresme proved already in the 14$^{\text{th}}$ century that the harmonic series does not converge. This simple proof is well known

$$\frac{1}{1} + (\frac{1}{2}) + (\frac{1}{3} + \frac{1}{4}) + (\frac{1}{5} + \frac{1}{6} + \frac{1}{7} + \frac{1}{8}) + ...$$

Every sum in parentheses amounts to at least 1/2. If we take infinitely many sums, we get infinitely many times 1/2 or more, so that the total sum is not finite. Counting the pairs of parentheses, we find that not less than $\aleph_0$ of them are necessary[7]. Counting the fractions we find that Nicole de Oresme used $2^{\aleph_0}$ natural numbers as denominators, not aware of Cantors celebrated theorem (1) according to which so many natural numbers are not available.

Several important definitions require more than $\aleph_0$ natural numbers. If, for instance, the sequence 1/n! cannot be calculated for infinitely many terms, i. e. up to $1/\aleph_0!$, then Euler's number e is not irrational let alone transcendental. So $\aleph_0!$ natural numbers must unavoidably exist to save the existence of e.

The first transcendental numbers, which J. Liouville found in 1844, are calculated by infinite series like $\tau = \Sigma 10^{-\nu!}$. Here $10^{\aleph_0!}$ natural numbers are needed. Would they not exist, there were no transcendental Liouville numbers and presumably no transcendental numbers at all. From this point of view, the existence of transcendental numbers proves the equivalence

$$|\mathbb{N}| = |\mathbb{R}| \tag{2}$$

which, if no transcendental numbers would exist, was also true, because in that case only the set of algebraic numbers with its cardinality $\aleph_0$ remained. Consequently, equivalence (2) is true in any case.

---

[7] A smaller, finite, number would not yield an infinite sum, because each one contributes less than 1. $\aleph_0$ terms are also required to establish the identities 0.999... = 1.000...  or  1/2 + 1/4 + 1/8 + ... = 1.

The asserted non-existence of a solution of

$$\lim_{y \to \infty} y = 2^x \quad \text{or} \quad \aleph_0 = 2^x \tag{3}$$

is questionable from the historical as well as from the mathematical point of view.

The historical point is the following. There have been always equations without solutions. One of them was $1 + x = 0$ until Leonardo of Pisa introduced negative numbers, another one was $1 + x^2 = 0$ and it remained so even longer, like negative logarithms and related notions. But it was never heard of, in any branch of mathematics, that an equation which once upon a time had had a solution, later on lost it. Now, in finite mathematics all equations of the form $y = 2^x$ have solutions. Before the advent of set theory, equation (3) had a solution too, namely $x = \infty$. Meanwhile there is no longer any solution at all. This is the only case in history where the set of solutions has shrunk.

The mathematical point is that some bijections have become impossible and that variables have to jump over large gaps. There are bijections which are well defined, or at least can be defined, over the whole sequence of natural numbers (so far it exists) like $n \leftrightarrow 2n$ and some others which are given in table 1, including their limit.

**Table 1.** Some everywhere defined bijections.

| $n$ | $2n$ | $n^2$ | $1/n$ |
|---|---|---|---|
| 1 | 2 | 1 | 1 |
| 2 | 4 | 4 | 1/2 |
| 3 | 6 | 9 | 1/3 |
| ... | ... | ... | ... |
| $\aleph_0$ | $\aleph_0$ | $\aleph_0$ | 0 |

Bijections like $n \leftrightarrow 2^n$ and some others given in table 2, however, are discontinuous.

$$2^n < \aleph_0 \quad \text{for} \quad n < \aleph_0 \quad \text{and} \quad 2^n = 2^{\aleph_0} > \aleph_0 \quad \text{for} \quad n = \aleph_0 \tag{4}$$

This bijection is undefined between $\aleph_0$ and $2^{\aleph_0}$. There is a gap where $2^n$ is neither finite nor infinite. But if $\aleph_0$ is considered a *number* which can be in trichotomy with natural numbers, then this discontinuity is highly suspect. $n \leftrightarrow 2^{2^{n!}}$ exhibits indefiniteness over even a larger gap. All the rows with question marks in table 2 show gaps of indefiniteness. They are unanswered by transfinite set theory, which, on the other hand, is ready to calculate hyper-inaccessible cardinals minutely. The common calculus simply would replace all question marks by $\infty$.

**Table 2.** Some bijections with gaps.

| $1/2^{n!}$ | $1/n!$ | $\log_2 n$ | $n$ | $2^n$ | $n!$ | $2^{n!}$ | $2^{2^{n!}}$ |
|---|---|---|---|---|---|---|---|
| 1/2 | 1 | 0 | 1 | 2 | 1 | 2 | 4 |
| 1/4 | 1/2 | 1 | 2 | 4 | 2 | 4 | 16 |
| 1/64 | 1/6 | $\log_2 3$ | 3 | 8 | 6 | 64 | $2^{64}$ |
| ... | ... | ... | ... | ... | ... | ... | ... |
| ? | ? | ? | ? | ? | ? | ? | $\geq \aleph_0$ |
| ... | ... | ... | ... | ... | ... | ... | ... |
| 0 | ? | ? | ? | ? | ? | $\geq \aleph_0$ | $\geq 2^{\aleph_0}$ |
| ... | ... | ... | ... | ... | ... | ... | ... |
| 0 | 0 | ? | ? | ? | $\geq \aleph_0$ | $\geq 2^{\aleph_0}$ | $\geq 2^{2^{\aleph_0}}$ |
| ... | ... | ... | ... | ... | ... | ... | ... |
| 0 | 0 | ? | ? | $\geq \aleph_0$ | $\geq \aleph_0$ | $\geq 2^{\aleph_0}$ | $\geq 2^{2^{\aleph_0}}$ |
| ... | ... | ... | ... | ... | ... | ... | ... |
| 0 | 0 | ? | $\geq \aleph_0$ | $\geq 2^{\aleph_0}$ | $\geq \aleph_0!$ | $\geq 2^{\aleph_0!}$ | $\geq 2^{2^{\aleph_0!}}$ |
| ... | ... | ... | ... | ... | ... | ... | ... |
| 0 | 0 | $\geq \aleph_0$ | $\geq 2^{\aleph_0}$ | $\geq 2^{2^{\aleph_0}}$ | $\geq 2^{\aleph_0}!$ | $\geq 2^{(2^{\aleph_0})!}$ | $\geq 2^{2^{(2^{\aleph_0})!}}$ |

Of course this table could easily be extended by loglog...log $n$ as well as by higher powers.

## 4. Enumerating the real numbers

In possession of at least $2^{\aleph_0}$ natural numbers and convinced of the non-existence of actually infinite numbers (of digits), a surjective and injective mapping of $\mathbb{N}$ on the real interval (0, 1] can be accomplished. Such a bijection is demonstrated in table 3 (others are possible).

**Table 3.** Enumerated list of all real numbers $r \in (0, 1]$ in binary representation.

| 1) | 0.1 | 2) | 0.01 | 4) | 0.001 | 8) | 0.0001 | ... |
|---|---|---|---|---|---|---|---|---|
| | | 3) | 0.11 | 5) | 0.101 | 9) | 0.1001 | ... |
| | | | | 6) | 0.011 | 10) | 0.0101 | ... |
| | | | | 7) | 0.111 | 11) | 0.1101 | ... |
| | | | | | | 12) | 0.0011 | ... |
| | | | | | | 13) | 0.1011 | ... |
| | | | | | | 14) | 0.0111 | ... |
| | | | | | | 15) | 0.1111 | ... |
| | | | | | | | | ... |

The numbers of table 3 are obtained by starting with 0.1 and then inserting 0 or 1 behind the decimal point of the numbers given in the preceding column. A subset filling one column of table 3 is always countable, because it surpasses the cardinality of the union of the preceding columns only by 1. Cantor calls it a simple and easy to prove theorem, that the union of a countable set of countable sets is countable.[8] Further table 3 contains all real numbers including all those with (potentially) infinite non-periodic strings of digits, each one appearing at a unique place of the sequence.[9] Therefore, it constitutes a correct bijection as in detail specified by Cantor[10] (though with respect to the set of algebraic numbers).

As there is only a countable set of intervals of length 1 in $\mathbb{R}$, the denumerability of the whole set $\mathbb{R}$ is obvious. Again we find equation (2) being true and inequality (1) being false.

---

[8] Hat man eine endliche oder abzählbar unendliche Menge von Mengen (E), (E'), (E''), ..., deren jede ihrerseits abzählbar ist, so ist auch die aus der Zusammenfassung aller Elemente von (E), (E'), (E''), ... hervorgehende Menge abzählbar. [1, p. 152] Da nämlich alle Bestandteile ... abzählbar sind und die Anzahl ... eine abzählbar unendliche ist, so folgt daraus die Abzählbarkeit ... [1, p. 160].

[9] Often at this point the question is raised which natural number is mapped on √2 or on π? This question is easily answered, once *all* bits of √2 or π are provided by the questioner.

[10] ... so daß zu jeder algebraischen Zahl ω eine bestimmte ganze positive Zahl ν und umgekehrt zu jeder positiven ganzen Zahl ν eine völlig bestimmte reelle algebraischen Zahl y gehört, daß also, um mit anderen Worten dasselbe zu bezeichnen, der Inbegriff (ω) in der Form einer unendlichen gesetzmäßigen Reihe $\omega_1, \omega_2, ... \omega_\nu,...$ gedacht werden kann, in welcher sämtliche Individuen von (ω) vorkommen und ein jedes von ihnen sich an einer bestimmten Stelle ..., welche durch den zugehörigen Index gegeben ist, befindet [1, p. 115].